\documentclass[12pt]{amsart}
\usepackage{geometry}                
\usepackage{graphicx}
\usepackage{pstricks,pst-poly,pst-node,pst-text,pst-plot,pstricks-add}

\usepackage{amssymb}
\newtheorem{prop}{Proposition}
\newcommand{\beq}{\begin{equation}}
\newcommand{\eeq}{\end{equation}}
\usepackage{pstricks,pst-poly,pst-node,pst-text,pst-plot,pstricks-add}

\title{Another Incarnation of the Lambert $W$ Function }
\author{Alexander Kheyfits}

\begin{document}
\maketitle
 
\begin{flushright} To the memory of Jonathan M. Borwein (1951-2016)
\end{flushright} 

\vspace{1cm}
\begin{center}
\today    \\
\end{center}
\vspace{.3cm}

\textbf{Abstract:} \emph{The Lambert $W$ function was introduced by Euler in 1779, but was not well-known until it was implemented in Maple, and the seminal paper \cite{CorGHJK} was published in 1996. In this note we describe a simple problem, which can be straightforwardly solved in terms of the $W$ function. }\\
\vspace{.5cm} 

At the recent seminar, Professor Bertram Kabak asked the following question: \\

\emph{The graphs of the exponential function $e^x$ and of its inverse, the natural logarithm $\ln x$ have no point in common, see Fig. 1. Consider an exponential function with a base $b,  \; 0<b,\, b\neq 1 $, 
\[y=b^x.\] 
For what values of the base $b$, have the graphs of this function and its inverse, $y(x)=\log_bx$ the common points, and how many? }  
\begin{figure}[hbt]
\psset{xunit=.5cm,yunit=.5cm}
\def\xlim{8.5}
\def\ylim{4.5}
\begin{center}
\begin{pspicture*}(-2.5,-2.5)(\xlim,\ylim)
\psaxes[Dx=2,Dy=2]{->}(0,0)(-2.5,-4.5)(\xlim,4.5)
\rput(8.2,.35){$x$}        \rput(.35,4.3){$y$}

\psplot[plotpoints=200]{-6}{1.5}{2.71828 x exp}
\psline[linestyle=dashed](-3,-3)(4,4)
\psplot[plotpoints=200]{.05}{8}{x ln }
\rput(1.8,3.5){$y=e^x$}        \rput(7,1.5){$y=\ln x$}    
\end{pspicture*}  
\end{center}
\caption{Graphs of exponential function $y=e^x$ and its inverse logarithmic function $y=\ln x$ never intersect. }
\end{figure}
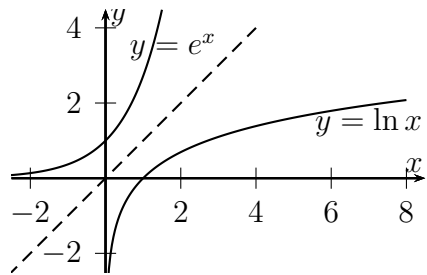

It turns out that this problem can be easily solved by making use of the Lambert $W$ function $W=w(z)$. The latter is given implicitly for every complex number $z$ by the equation   
\beq we^w=z, \eeq
where we follow the standard-now notation $w=w(z)$ \cite{CorGHJK}.
 
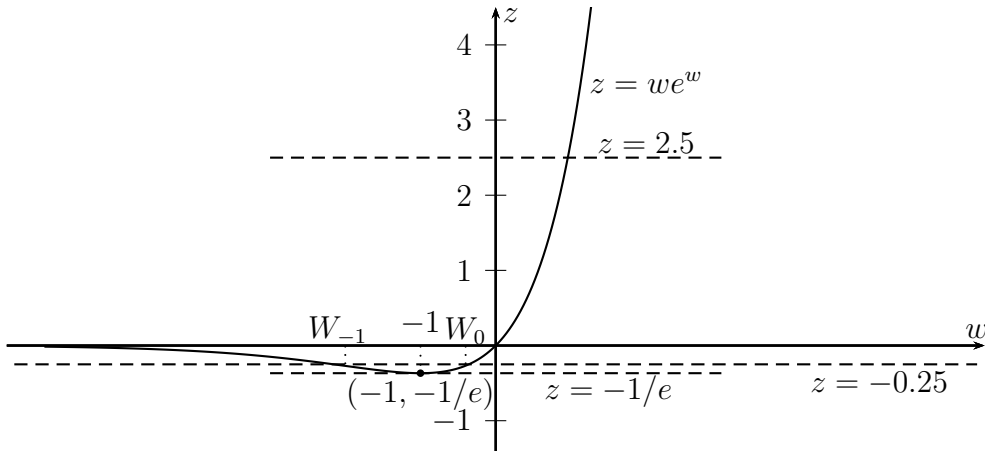
\begin{figure}[hbt]
\psset{xunit=1cm,yunit=1cm}
\def\xlim{6.5}
\def\ylim{4.5}
\begin{center}
\begin{pspicture*}(-\xlim,-1.4)(\xlim,\ylim)
\psaxes[Dx=14,Dy=1]{->}(0,0)(-\xlim,-4.5)(\xlim,4.5)
\rput(6.4,.2){$w$}        \rput(.2,4.4){$z$}
\psplot[plotpoints=200]{-6}{1.5}{x 2.71828 x exp mul}
\psline[linestyle=dashed](-3,2.5)(3,2.5)
\psline[linestyle=dashed](-6.4,-.25)(6.4,-.25)
\psline[linestyle=dashed](-3,-0.367879)(3,-.367879)
\psline[linestyle=dotted](-1,0)(-1,-.39)
\psline[linestyle=dotted](-.4,0)(-.4,-.25)
\psline[linestyle=dotted](-2,0)(-2,-.25)
\rput(2,3.5){$z=we^w$}        \rput(-1,.3){$-1$}        \rput(1.5,-.6){$z=-1/e$}     \rput(2,2.7){$z=2.5$}    \rput(5.1,-.5){$z=-0.25$}       \rput(-.4,.2){$W_0$}                       \rput(-2.1,.2){$W_{-1}$}   
\pscircle*(-1,-.367879){.05}   \rput(-1,-.65){$(-1,-1/e)$}  
\end{pspicture*}   
\end{center}
\caption{Function $z=we^w$. The dashed lines are the horizontal tangent line $z=-1/e$ and the horizontal lines $z=2.5$ and $z=-0.25$. }
\end{figure}

After the $W$ was implemented in Maple and other major CASs, and the article \cite{CorGHJK} has been published, the number of publications devoted to the function, has grown dramatically\footnote{On Jan. 8, 2017, the query 'Lambert $W$ function' returned 302 references on MathSciNet, about 450,000 references on Google, and more than 1,200,000 references on Google Scholar.}. The function enjoys so many applications, that together with the logarithmic function, the $W$ should be in a toolbox of any researcher. Some authors claim that the $W$ is an elementary function. Leave it to the individual judgment, whether and in what sense the Lambert $W$ function is elementary, we show how the $W$ naturally occurs in an elementary problem.  \\

The Lambert $W$ is a many-valued analytic function, therefore, its complete study and, in particular, description of its single-valued branches, should be done in complex-analytic terms, which is beyond the scope of this brief note, see, e.g., \cite{CorGHJK} or \cite{Kh}, where the closed-form representation of all the branches of $W$ in terms of contour integrals was derived.   

We mention only few necessary properties of the $W$ function, referring the reader to the papers above. Begin by graphing the left-hand side of (1), see Fig. 2. The function $f(w)=w e^w$ has the global minimum $-1/e$ at $w=-1$. For every $z\geq 0$, equation $z=we^w$ has exactly one real root, that is, the Lambert $W$ function has one real-valued branch. If $-1/e<z<0$, equation (1) has two real roots, $-1<W_0<0$ and $W_{-1}<-1$. When $z=-1/e$, these two roots merge in the double root $w=-1$. When $z<-1/e$, these roots disappear from the graph, because they became complex numbers. 

The function $f(w)=w e^w - z$ is an entire function for every complex $z$; such functions are called quasi-polynomials. It is distinct from $e^w$, therefore, for any complex $z\neq 0$ this function has infinitely many roots. Since every root generates its own single-valued branch , the inverse function, Lambert $W$ has infinitely many branches, and only two of them are real-valued on the real axis. They are conventionally called the principal branch, $W_0(z)$, which is real  for $z\geq -1/e$, and the branch $W_{-1}(z)$, real-valued for $-1/e\leq z <0$.  

The graph of $W_0(z)$ is the mirror reflection, with respect to the bisectrix $w=z$, of the right half,
$w\geq -1,\; z\geq -1/e$, of the graph in Fig. 2. The graph of $W_{-1}(z)$ is the mirror reflection of the left half of the same graph, $w\leq -1,\; -1/e \leq z \leq 0$, with respect to the bisectrix $z=w$. 

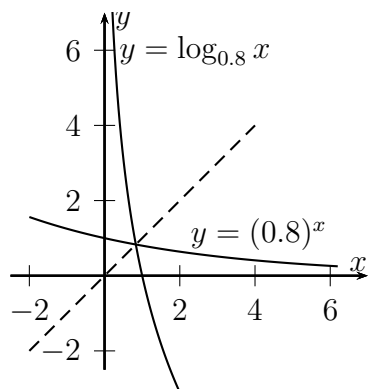
\begin{figure}[hbt]
\psset{xunit=.5cm,yunit=.5cm}
\def\xlim{7}
\def\ylim{7}
\begin{center}
\begin{pspicture*}(-2.5,-3)(\xlim,\ylim)
\psaxes[Dx=2,Dy=2]{->}(0,0)(-2.5,-2.5)(\xlim,\ylim)
\rput(6.75,.3){$x$}        \rput(.5,6.8){$y$}
%\rput(-.4,.4){$\mathcal{O}$}

\psplot[plotpoints=200]{-2}{6.2}{.8 x exp}
\psline[linestyle=dashed](-2,-2)(4,4)
\rput(2.4,6){$y=\log_{0.8} x$}    
\rput(4.1,1.1){$y=(0.8)^x$}          
\psplot[plotpoints=200]{.2}{7}{x ln .8 ln div}
\end{pspicture*}  
\end{center}
\caption{The functions $y=(0.8)^x$ and $y=\log_{0.8} x$ have the unique intersection point. The dashed line is the bisectrix $y=x$.}
\end{figure}

Now we can take up the question above. First, let be $0<b<1$. The graphs of a function and its inverse are symmetrical with respect to the bisectrix $z=w$, in the notations of Fig. 2, or the bisectrix $y=x$ in the notations of Fig. 3. Therefore, the intersection points, if any, must belong to the bisectrix, that is, satisfy the equation $b^w=w$, or $(-\ln b)we^{-(\ln b)w}=-\ln b$. The right hand side here is positive because $b<1$, hence from Fig. 2 it follows that for any $0<b<1$ the equation has exactly one solution, or there exists the unique intersection point of $b^w$ and $\log_bw$. It is worth repeating that the intersection point is given by the principal branch $W_0(\ln (1/b)$ of the Lambert $W$ function. For instance, the case $b=0.8$ is shown in Fig. 3, where the exponential and logarithmic curves intersect at the point 
$W_0(-\ln 0.8)/(-\ln 0.8)\approx 0.83$. If we depart from the equation $w=\log_b w$, we arrive at the same conclusion. 

The case $b>1$ is a bit more interesting. The same reasoning leads to the equation $x=b^x$, or to equation (1), 
\[we^w=z, \]
where we are to set $w=-x\ln b$ and $z=- \ln b <0$. Since now $b>1$ and $x>0$, then $w=-x\ln b<0$. Therefore, if $-1/e<z=-\ln b <0$, that is, if $1<b<e^{1/e}$, then there are two intersections, $W_0(z)$ and $W_{-1}(z)$, given by the principal branch $W_0$ and the preceding branch $W_{-1}$ of the $W$ function, see Fig. 4.

If $-1/e=z$, where $z=-\ln b$, that is, $b=e^{1/e}$, then the two points merge into the point of tangency, see Fig. 5. We can state the conclusion as follows.  \\

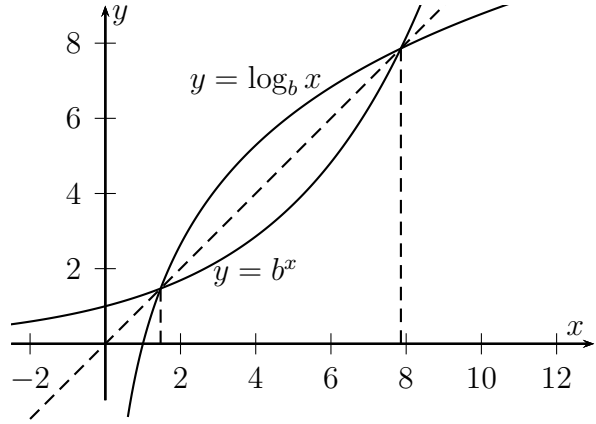
\begin{figure}[hbt]
\psset{xunit=.5cm,yunit=.5cm}
\def\xlim{13}
\def\ylim{9}
\begin{center}
\begin{pspicture*}(-2.5,-3)(\xlim,\ylim)
\psaxes[Dx=2,Dy=2]{->}(0,0)(-2.5,-1.5)(\xlim,\ylim)
\rput(12.5,.4){$x$}        \rput(.4,8.75){$y$}

\psplot[plotpoints=200]{-3}{14}{1.3 x exp}
\psline[linestyle=dashed](-2,-2)(12,12)
\psline[linestyle=dashed](7.86,0)(7.86,7.86)
\psline[linestyle=dashed](1.471,0)(1.471,1.471)
\rput(4,7){$y=\log_b x$}    
\rput(4,1.9){$y=b^x$}          
\psplot[plotpoints=200]{.6}{11}{x ln 1.3 ln div}
\end{pspicture*}  
\end{center}
\caption{The functions $y=b^x$ and $y=\log_b x$ with $b=1.3$ have two points of intersection.}
\end{figure}

\begin{prop} For $0<b<1$, the graph of the exponential function $y=b^x$ and the graph of its inverse $y=\log_b x $ have the unique point of intersection, given by the principal branch of the Lambert $W$ function $W_0(-\ln b)/(-\ln b)$. The graphs of the exponential function $y=b^x$ and of its inverse 
$y=\log_bx $ have a point of tangency if and only if $b=e^{1/e}$. The tangency point is $x_t$ with 
$-x_t \ln b =-1$, that is, $x_t=e$. These graphs intersect exactly twice if and only if $1<b<e^{1/e}$. The two points of intersection are given by the two branches, $W_0(-1/\ln b)$ and $W_{-1}(-1/\ln b)$. For example, if $b=1.3$, the intersection points are (Fig.4)
\[-W_0(-1/\ln 1.3)/\ln 1.3 \approx 0.386/\ln 1.3 \approx 1.47 \]
and 
\[-W_{-1}(-1/\ln 1.3)\approx 2.061 / \ln 1.3 \approx 7.86. \]  

For $b> e^{1/e}$, these graphs do not intersect (Fig. 1).         $\hfill \qed$
\end{prop}

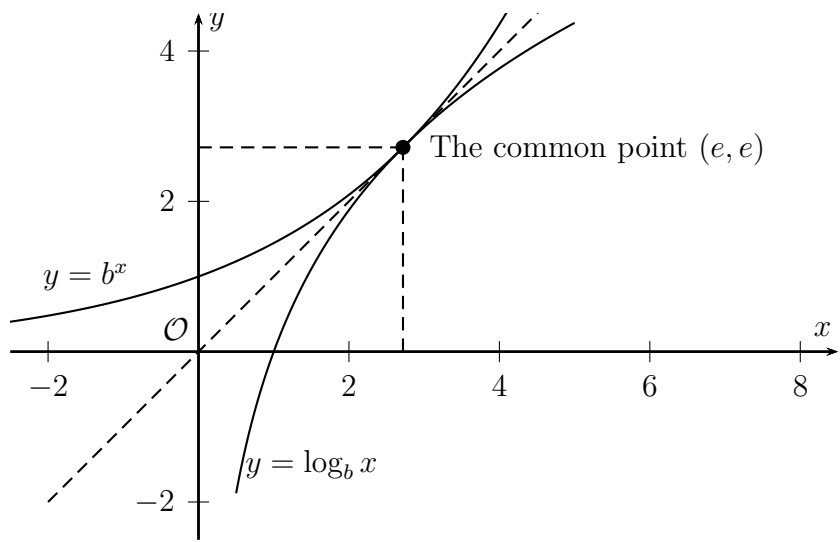
\begin{figure}[hbt]
\psset{xunit=1cm,yunit=1cm}
\def\xlim{8.5}
\def\ylim{4.5}
\begin{center}
\begin{pspicture*}(-2.5,-2.5)(\xlim,\ylim)
\psaxes[Dx=2,Dy=2]{->}(0,0)(-2.5,-2.5)(\xlim,4.5)
\rput(8.3,.3){$x$}        \rput(.25,4.4){$y$}
\rput(-.3,.3){$\mathcal{O}$}
\rput(2.71828,2.71828){$\pscircle*{.1}$}
\psplot[plotpoints=200]{-6}{6}{2.71828 x 2.71828 div exp}
\psline[linestyle=dashed](-2,-2)(5,5)
\rput(1.5,-1.5){$y=\log_b x$}    
\rput(-1.5,1){$y=b^x$}          \rput(5.3,2.7){The common point $(e,e)$}   
\psplot[plotpoints=200]{.5}{5}{x ln 2.71828 mul}
\psline[linestyle=dashed](0,2.71828)(2.71828,2.71828)
\psline[linestyle=dashed](2.71828,0)(2.71828,2.71828)
\end{pspicture*}  
\end{center}
\caption{The functions $y=b^x$ and $y=\log_b x\;$ with $b=e^{1/e} \approx 1.44$ have a point of tangency.}
\end{figure}

The discussion leads to the following natural problem. \\

{\bf Problem.} \emph{Describe those continuous increasing functions $f(x),\; 0<x<\infty$, whose graphs have a) one, b) two, c) none, d) $n\geq 3$ points of intersection with the inverse function $f^{-1}$. The additional assumption of convexity will, definitely, simplify the problem and likely, lead to the functions $e^{f(x)}$.} \\

{\bf Acknowledgement.} The author is thankful to Professor German Kalugin for useful remarks, which help to improve the manuscript.

\end{document}